\documentclass[11pt]{article}
\usepackage{mathrsfs}
\usepackage{amsfonts}
\usepackage{amsmath,amssymb}
\begin{document}
\title{\bf Fixed points of commutative L\"{u}ders operations}
\author {{Liu Weihua,\,\,  Wu Junde\thanks{Corresponding author: E-mail: wjd@zju.edu.cn}}\\
{\small\it Department of Mathematics, Zhejiang University, Hangzhou
310027, P. R. China}}

\date{}
\maketitle

{\small\it {\bf Abstract.} This paper verifies a conjecture posed in
a pair of papers on the fixed point sets for a class of quantum
operations. Specifically, it is proved that if a quantum operation
has mutually commuting operation elements that are effects forming a
resolution of the identity, then the fixed points set of the quantum
operation is exactly the commutant of the operation elements.}

\vskip 0.2 in


\vskip 0.2 in

\noindent {\bf 1. Introduction}

\vskip 0.2 in

\noindent Let $H$ be a complex Hilbert space, ${\cal B}(H)$ be the
bounded linear operator set on $H$. If $A\in {\cal B}(H)$ and $0\leq
A\leq I$, then $A$ is called a {\it quantum effect} on $H$. Each
quantum effect can be used to represent a yes-no measurement that
may be unsharp ([1-6]). The set of all quantum effects on $H$ is
denoted by ${\cal E}(H)$, the set of all orthogonal projection
operators on $H$ is denoted by ${\cal P}(H)$. Each element $P$ of
${\cal P}(H)$ can be used to represent a yes-no measurement that is
sharp ([1-6]). Let ${\cal T}(H)$ be the set of all trace class
operators on $H$ and ${\cal D}(H)$ the set of all density operators
on $H$, i.e., ${\cal D}(H)=\{\rho: \rho\in {\cal T}(H), \rho\geq 0,
tr (S)=1\}$. Each element $\rho$ of ${\cal D}(H)$ represents a state
of the quantum system $H$.

\vskip 0.1 in

Let ${\cal A}=\{E_i\}_{i=1}^n\subseteq {\cal E}(H)$ be a quantum
measurement, that is $\sum_{i=1}^nE_i^2=I$ in the strong operator
topology, where $1\leq n\leq\infty$, then the probability of outcome
$E_i$ is measured in the state $\rho$ is given by $tr(\rho E_{i})$,
and the new quantum state after the measurement ${\cal A}$ is
performed is defined by
$$\Phi(\rho)=\sum\limits_{i=1}^{n}E_{i}\rho E_{i}.$$

\vskip 0.1 in

Note that $\Phi: \rho\rightarrow \sum\limits_{i=1}^{n}E_{i}\rho
E_{i}$ defined a transformation on the state set ${\cal D}(H)$, we
call it the {\it L\"{u}ders transformation} ([6-7]). In physics, the
question whether a state $\rho$ is not disturbed by the measurement
${\cal A}=\{E_i\}_{i=1}^n$ becomes equivalent to the fact that
$\rho$ is a solution of the equation
$$\Phi(\rho)=\sum\limits_{i=1}^{n}E_{i}\rho
E_{i}=\rho.$$

It was showed in [8] that the measurement ${\cal A}=\{E_i\}_{i=1}^2$
does not disturb $\rho$ if and only if $\rho$ commutes with each
$E_i$, $i=1,2$.

\vskip 0.1 in

Moreover, if we define the {\it L\"{u}ders quantum operation
$\Phi_{{\cal A}}$} on ${\cal B}(H)$ as following:
$$\Phi_{{\cal A}}: {\cal B}(H)\rightarrow {\cal B}(H), \,\, B\rightarrow\Phi_{{\cal A}}(B)=\sum\limits_{i=1}^{n}E_{i}B
E_{i},$$ then an interesting problem is that if $B\in {\cal B}(H)$
is a fixed point of $\Phi_{{\cal A}}$, that is, $\Phi_{{\cal
A}}(B)=\sum\limits_{i=1}^{n}E_{i}B E_{i}=B,$ then $B$ commutes with
each $E_i$ ? $i=1, 2, \cdots, n$.

\vskip 0.1 in

In [9-10], we knew the conclusion is true if $H$ is a finite
dimensional complex Hilbert space. In [9-11], it was showed that the
conclusion is not true when $n=5$ or $n=3$ for infinite dimensional
complex Hilbert space. Thus, the general conclusion for infinite
dimensional cases is false. On the other hand, Busch and Singh in
[8] showed that for $n=2$ the conclusion is true for all complex
Hilbert spaces. Note that in this case, $E_1^2+E_2^2=I$, so
$E_1E_2=E_2E_1$, that is, ${\cal A}=\{E_1, E_2\}$ is commutative.
This motivated Arias, Gheonda, Gudder and Nagy to conjecture when
${\cal A}=\{E_i\}_{i=1}^n\subseteq {\cal E}(H)$ is commutative, then
the conclusion is true, that is, the fixed point set of $\Phi_{{\cal
A}}$ is exactly the commutant $\mathcal{A}'$ of the operation
elements $\mathcal{A}=\{E_{i}\}_{i=1}^n$. Moreover, Nagy in [12]
showed that if the conjecture is true, then
$$\Phi_{{\cal A}}(E)=\sum\limits_{i=1}^{n}E_{i}E
E_{i}=I-E$$ has the unique solution $\frac{1}{2}I$ in
${\cal E}(H)$, in physics, it showed that if the measurement ${\cal
A}$ disturbs the quantum effect $E$ completely into its supplement
$I-E$, then $E$ has to be $\frac{1}{2}I$.

\vskip 0.1 in

As showed in [13-16], the structures of fixed point sets of quantum
operations have important applications in quantum information
theory, in particular, in [15, Theorem 3], the fixed point set is a
matrix algebra which share an elegant structure, played a central
role in identifying the protected structures.

\vskip 0.1 in

In this paper, by using the spectral theory of self-adjoint
operators, we prove the conjecture affirmatively. Moreover, when
${\cal A}=\{E_i\}_{i=1}^{n}\subseteq {\cal E}(H)$ is commutative and
$F=\sum_{i=1}^{n}E_i^2<I$, we also obtain a nice conclusion. Note
that the von Neumann algebra $\mathcal{N}$ generated by
$\{E_i\}_{i=1...n}$ is Abelian which can be embed into a maximal
Abelian von Neumann algebra. Since a maximal Abelian von Neumann
algebra $\mathcal{M}$ on a separable Hilbert space is always a
direct sum of $\mathcal{M}_1$ and $\mathcal{M}_2$. Here
$\mathcal{M}_1$ is isometric to $\bigoplus\limits_{i=1}^{\infty}
C_i$ and $\mathcal{M}_2$ is isometric to $L_{\infty}(B)$, where $B$
is a compact subset of the real number set $R$. Thus, $\mathcal{A}'$
has the form $\bigoplus\limits_{i=1}^{\infty} M_k\otimes
1_{n_k}\bigoplus L_{\infty}(C) $, where $C$ is a subset of $B$ and
$M_k$ is a matrix algebra whose dimension is k and $n_k$ ranges from
0 to $\infty$ ([17]). So our conclusions is analogous with the
finite dimensional cases' concise shape in Theorem 3 in [15].

\vskip 0.2 in

\noindent {\bf 2. Element lemmas and proofs}

\vskip 0.2 in

Let $1\leq n<\infty$ and $\mathcal{A}=\{E_i\}_{i=1}^{n}\subseteq
{\cal E}(H)$ be commutative. Firstly, for each $E_i, 1\leq i\leq n$,
we have the spectral representation theorem:
$$E_i=\int\limits_0^1 \lambda dF^{(i)}_\lambda,$$
where $\{F^{(i)}_\lambda\}_{\lambda\in \mathbb{R}}$ is the identity
resolution of $E_i$ satisfying that $\{F^{(i)}_\lambda\}_{\lambda\in
\mathbb{R}}$ is right continuous in the strong operator topology and
$F^{(i)}_\lambda=0$ if $\lambda< 0$ and $F^{(i)}_\lambda=I$ if
$\lambda\geq ||E_i||$, moreover, for each $\lambda\in \mathbb{R}$,
$F^{(i)}_{\lambda}=P^{E_i}(-\infty, \lambda]$, where $P^{E_i}$ is
the spectral measure of $E_i$ ([17]). Now, for fixed integers $m,
k_1, k_2, \cdots, k_n$, we denote
$$F^m_{k_1,\cdots,k_n}=P^{E_1}(\frac{k_1}{m},\frac{k_1+1}{m}]\cdots
P^{E_n}(\frac{k_n}{m},\frac{k_n+1}{m}].$$ Since $E_i$ and $E_j$ are
commutative for any $i,j=1, 2, \cdots, n$, so $F^m_{k_1...k_n}$ is a
well-defined orthogonal projection operator.

\vskip 0.1 in

{\bf Lemma 2.1.} Let $1\leq n<\infty$, ${\cal
A}=\{E_i\}_{i=1}^{n}\subseteq {\cal E}(H)$ be commutative and $B\in
{\cal B}(H)$. If for any integers $m$ and $k_1, k_2, \cdots,k_n$,
$B$ commutes with $F^m_{k_1,\cdots,k_n}$, then $B$ is commutative
with each $E_i$ in $\mathcal {A}=\{E_i\}_{i=1}^{n}$.

{\bf Proof}. For each rational number $q=\frac{p}{l}$, where $p,l$
are integers. If $\frac{p}{l}<0$, then $F^{(i)}_{\frac{p}{l}}=0$, if
$\frac{p}{l}\geq 1$, then $F^{(i)}_{\frac{p}{l}}=I$. Let $l>p\geq0$,
so $0\leq \frac{p}{l}<1$. Then
$F^{(i)}_{\frac{p}{l}}=P^{E_i}(\frac{-1}{l}, 0]+P^{E_i}(0,
\frac{1}{l}]+\cdots + P^{E_i}(\frac{p-1}{l}, \frac{p}{l}]$, thus, we
can prove easily that
$$F^{(i)}_{\frac{p}{l}}=\sum\limits_{k_i<p}(\sum\limits_{k_1,\cdots,k_{i-1},k_i,k_{i+1},\cdots,k_n}F^l_{k_1,\cdots,k_n}).$$
So, for each rational number $q=\frac{p}{l}$,
$F^{(i)}_{\frac{p}{l}}$ commutes with $B$, note that
$\{F^{(i)}\}_{\lambda\in \mathbb{R}}$ is right continuous in the
strong operator topology, so $B$ commutes with each $E_i$, $i=1, 2,
\cdots, n$.

\vskip 0.1 in

{\bf Lemma 2.2.} Let $1\leq n<\infty$, ${\cal
A}=\{E_i\}_{i=1}^{n}\subseteq {\cal E}(H)$ be commutative and $B\in
{\cal B}(H)$. If $B$ does not commute with some $E_{i_0}$ in
$\mathcal {A}$, then there are integers $m$, $k_1, k_2, \cdots, k_n$
and $k_1', k_2', \cdots, k_n'$, such that $k_i\neq k_i'$ for at
least one $i$ and $F^m_{k_1, k_2, \cdots, k_n}BF^m_{k'_1, k_2',
\cdots, k'_n}\neq 0$.

{\bf Proof.} Without of losing generality, we suppose that $B$ does
not commute with $E_1$. By Lemma 2.1, there are integers $m$ and
$k_1, k_2, \cdots, k_n$ such that $F^m_{k_1, k_2, \cdots, k_n}B\neq
F^m_{k_1, k_2, \cdots, k_n}BF^m_{k_1, k_2, \cdots, k_n}$ or
$BF^m_{k_1, k_2, \cdots, k_n}\neq F^m_{k_1, k_2, \cdots,
k_n}BF^m_{k_1, k_2, \cdots, k_n}$. If $F^m_{k_1, k_2, \cdots,
k_n}B\neq F^m_{k_1, k_2, \cdots,k_n}BF^m_{k_1, k_2, \cdots, k_n}$,
then there exists integers $k_1', k_2', \cdots, k_n'$, $k_i\neq
k_i'$ for at least one $i$ such that $F^m_{k_1, k_2, \cdots,
k_n}BF^m_{k'_1, k_2', \cdots, k'_n}\neq 0$. In fact, if not, we will
get that $$F^m_{k_1, k_2, \cdots, k_n}B=\sum\limits_{k_1', k_2',
\cdots, k'_n}F^m_{k_1, k_2, \cdots, k_n}BF^m_{k'_1, k_2', \cdots,
k'_n}= F^m_{k_1, k_2, \cdots, k_n}BF^m_{k_1, k_2, \cdots, k_n}.$$
This is a contradiction. Similarly, if $BF^m_{k_1, k_2, \cdots,
k_n}\neq F^m_{k_1, k_2, \cdots, k_n}BF^m_{k_1, k_2, \cdots, k_n}$,
we will also get the same conclusion. The lemma is proven.

Moreover, we have a stronger conclusion in the following.

\vskip 0.1 in

{\bf Lemma 2.3.} Let $A\in {\cal E}(H)$ and $B\in {\cal B}(H)$. If
$B$ does not commute with $A$, then there exists integer $m,k,j$
with $|k-j|\geq 2$ such that
$$P^{A}(\frac{k}{m},\frac{k+1}{m}]BP^{A}(\frac{j}{m},\frac{j+1}{m}]\neq
0.$$

{\bf Proof.} By lemma 2.2, we can find $k_1\neq j_1$ such that
$C=P^{A}(\frac{k_1}{m},\frac{k_1+1}{m}]BP^{A}(\frac{j_1}{m},\frac{j_1+1}{m}]\neq
0$. If $|k_1-j_1|\geq 2$, then we get the $m,k,j$ satisfy the lemma.
If $j_1=k_1+1$, we replace $m$ by $2m$ and let $k_2=2k_1$,
$j_2=2j_1$. Then
$$P^{A}(\frac{k_1}{m},\frac{k_1+1}{m}]=
P^{A}(\frac{k_2}{2m},\frac{k_2+1}{2m}]+P^{A}(\frac{k_2+1}{2m},\frac{k_2+2}{2m}],$$
$$P^{A}(\frac{j_1}{m},\frac{j_1+1}{m}]=
P^{A}(\frac{j_2}{2m},\frac{j_2+1}{2m}]+P^{A}(\frac{j_2+1}{2m},\frac{j_2+2}{2m}].$$

Now we consider $k_2,k_2+1$ and $j_2,j_2+1$, if we still can not
take $|k-j|\geq 2$ satisfy the conclusion, then
$$P^{A}(\frac{k_2}{2m},\frac{k_2+1}{2m}]BP^{A}(\frac{j_2}{2m},\frac{j_2+1}{2m}]=0,$$
$$P^{A}(\frac{k_2}{2m},\frac{k_2+1}{2m}]BP^{A}(\frac{j_2+1}{2m},\frac{j_2+2}{2m}]=0,$$
$$P^{A}(\frac{k_2+1}{2m},\frac{k_2+2}{2m}]BP^{A}(\frac{j_2+1}{2m},\frac{j_2+2}{2m}]=0.$$
So we have
$C=P^{A}(\frac{k_2+1}{2m},\frac{k_2+2}{2m}]BP^{A}(\frac{j_2}{2m},\frac{j_2+1}{2m}]$.

Keep on this way, then we will find the integers $k,j$ which satisfy
the conclusion or we get a sequence $\{p_i
,p_i+1,2^{i-1}m\}_{i=1}^{\infty}$ such that $p_i+1=2^{i-1}j_1$ and
$C=P^{A}(\frac{p_i}{2^{i-1}m},\frac{p_i+1}{2^{i-1}m}]BP^{A}(\frac{p_i+1}{2^{i-1}m},\frac{p_i+2}{2^{i-1}m}]$.
If the first case happens, then we proved the lemma. If the second
case happens, note that
$$\bigcap\limits_{i=1}^{\infty}(\frac{p_i+1}{2^{i-1}m},\frac{p_i+2}{2^{i-1}m}]=\emptyset,$$
and
$$\bigcap\limits_{i=1}^{\infty}(\frac{p_i}{2^{i-1}m},\frac{p_i+1}{2^{i-1}m}]=\{\frac{j_1}{m}\},$$ so
$\lim\limits_{i\rightarrow
\infty}P^{A}(\frac{p_i}{2^{i-1}m},\frac{p_i+1}{2^{i-1}m}]=P^{A}\{\frac{j_1}{m}\}$
 and $\lim\limits_{i\rightarrow
\infty}P^{A}(\frac{p_i+1}{2^{i-1}m},\frac{p_i+2}{2^{i-1}m}]=0$ in
strong operator topology, thus,
$$\lim\limits_{i\rightarrow \infty}P^{A}(\frac{p_i}{2^{i-1}m},\frac{p_i+1}{2^{i-1}m}]BP^{A}(\frac{p_i+1}{2^{i-1}m},\frac{p_i+2}{2^{i-1}m}]=0$$
in strong operator topology ([17]). But for each positive integer
$i$,
$$C=P^{A}(\frac{p_i}{2^{i-1}m},\frac{p_i+1}{2^{i-1}m}]BP^{A}(\frac{p_i+1}{2^{i-1}m},\frac{p_i+2}{2^{i-1}m}],$$
so we get that $C=0$, this is a contradiction and the lemma is
proved in this case.

If $k_1+1=j_1$, we just need to take all the above calculations in
adjoint and interchange the index j and k. The proof is similar,
thus, we proved the lemma.

\vskip 0.1 in

{\bf Lemma 2.4.} Let $1\leq n < \infty$, ${\cal
A}=\{E_i\}_{i=1}^{n}\subseteq {\cal E}(H)$ be commutative and
$\sum_{i=1}^{n}E_i^2\leq I$. If $X\in B(H)$ is not commutative with
$E_1$, then there exists positive integer $m$ such that for each
positive integer $p$, there exist projection operators $P,Q\in
\mathcal{A}'$, $PQ=0$, $Y=PXQ\not=0$, and
$$\frac{\|Y\|-\|\Phi_{\cal A}(Y)\|}{\|Y\|}\geq
\frac{p^2-4\sqrt{n}mp-2n}{2(pm)^2}.$$

{\bf Proof.} Since $X$ does not commute with $E_1$, it follows from
Lemma 2.3 that there exists integers $m,k,j$ such that $|k-j|\geq 2$
and
$P^{E_1}(\frac{k}{m},\frac{k+1}{m}]XP^{E_1}(\frac{j}{m},\frac{j+1}{m}]\neq
0$. Note that
$$P^{E_1}(\frac{k}{m},\frac{k+1}{m}]XP^{E_1}(\frac{j}{m},\frac{j+1}{m}]
=\sum\limits_{k_{2},\cdots,k_n}\sum\limits_{k'_{2},\cdots,k'_n}F^m_{k,k_2,\cdots,k_n}XF^m_{j,k_2'\cdots,k_n'},$$
so there exist $k, k_2,\cdots, k_n$ and $j, k_2', \cdots, k_n'$ such
that $|k-j|\geq2$ and $$F^m_{k, k_2, \cdots, k_n}XF^m_{j, k_2',
\cdots, k_n'}\neq 0.$$ Let $P_0=F^m_{k, k_2, \cdots, k_n}$,
$Q_0=F^m_{j, k_2', \cdots, k_n'}$, $Y_0=P_0XQ_0$. Then $P_0$ and
$Q_0$ are projection operators and $P_0,Q_0\in \mathcal{A}'$,
$P_0Q_0=0$, $Y_0=P_0XQ_0\neq 0$. Moreover, for each $i=1, 2, \cdots,
n$, if we denote $k_1=k$, $k'_1=j$, then
\begin{equation}
\begin{array}{rcl}
 \|E_iY_0E_i\|&=&\|E_iP^{Ei}(\frac{k_i}{m},\frac{k_i+1}{m}]Y_0P^{E_i}(\frac{k_i'}{m},\frac{k_i'+1}{m}]E_i\|\\
  &\leq&\|E_iP^{E_i}(\frac{k_i}{m},\frac{k_i+1}{m}]\|\|Y_0\|\|P^{E_i}(\frac{k_i'}{m},\frac{k_i'+1}{m}]E_i\|\\
  &\leq&\frac{k_i+1}{m}\|Y_0\|\frac{k_i'+1}{m}\\
  &=&\frac{k_i+1}{m}\frac{k_i'+1}{m}\|Y_0\|.\\
\end{array}
\end{equation} Thus, we have
\begin{equation}
\|\sum\limits_{i=1}^{n}E_iY_0E_i\|
\leq\sum\limits_{i=1}^{n}\|E_iY_0E_i\|
\leq(\sum\limits_{i=1}^{n}\frac{k_ik'_i}{m^2}+\sum\limits_{i=1}^{n}\frac{k_i+k'_i}{m^2}+\frac{n}{m^2})\|Y_0\|.\\
\end{equation}
Since $\sum\limits_{i=1}^{n}E_i^2\leq I$ and
\begin{equation}
\begin{array}{rcl}
F^m_{k,k_2,\cdots,k_n}(I-\sum\limits_{i=1}^{n}E_i^2)
         &=&F^m_{k,k_2,\cdots,k_n}-F^m_{k,k_2,\cdots,k_n}\sum\limits_{i=1}^{n}E_i^2\\
         &\leq& F^m_{k,k_2,\cdots,k_n}-\sum\limits_{i=1}^{n}\frac{k_i^2}{m^2}F^m_{k,k_2,\cdots,k_n}\\
         &=&
         (1-\sum\limits_{i=1}^{n}\frac{k_i^2}{m^2})F^m_{k,k_2,\cdots,k_n},
\end{array}
\end{equation} so, we have $\sum\limits_{i=1}^{n}k_i^2\leq m^2$.
Similarly, we have also $\sum\limits_{i=1}^{n}k_i'^2\leq m^2.$
Moreover, note that
\begin{equation}
\begin{array}{rcl}
2m^2(1-\sum\limits_{i=1}^{n}\frac{k_ik'_i}{m^2}-\sum\limits_{i=1}^{n}\frac{k_i+k'_i}{m^2}-\frac{n}{m^2})
&=&m^2+m^2-2\sum\limits_{i=1}^{n}k_ik'_i-2\sum\limits_{i=1}^{n}(k_i+k'_i)-2n\\
&\geq&\sum\limits_{i=1}^{n}k_i^2+\sum\limits_{i=1}^{n}k_i'^{2}-2\sum\limits_{i=1}^{n}k_ik'_i-2\sum\limits_{i=1}^{n}(k_i+k'_i)-2n\\
&=&\sum\limits_{i=1}^{n}(k_i-k'_i)^2-2\sum\limits_{i=1}^{n}(k_i+k'_i)-2n\\
&\geq&(k_1-k'_1)^2-2\sum\limits_{i=1}^{n}(k_i+k'_i)-2n,
\end{array}
\end{equation} and $(\sum\limits_{i=1}^{n}k_i)^2\leq
n(\sum\limits_{i=1}^{n}k_i^2)\leq nm^2$,
$(\sum\limits_{i=1}^{n}k_i')^2\leq
n(\sum\limits_{i=1}^{n}k_i'^2)\leq nm^2$, we have
\begin{equation}
2m^2(1-\sum\limits_{i=1}^{n}\frac{k_ik'_i}{m^2}-\sum\limits_{i=1}^{n}\frac{k_i+k'_i}{m^2}-\frac{n}{m^2})
\geq (j-k)^2-4\sqrt{n}m-2n.\\
\end{equation}

On the other hand, it follows from
$$\|Y_0\|-\|\sum\limits_{i=1}^{n}E_iY_0E_i\|\geq
\|Y_0\|-\sum\limits_{i=1}^{n}\|E_iY_0E_i\|$$$$\geq
[1-(\sum\limits_{i=1}^{n}\frac{k_ik'_i}{m^2}+\sum\limits_{i=1}^{n}\frac{k_i+k'_i}{m^2}+\frac{n}{m^2})]\|Y_0\|$$
and (5) that
$$\frac{\|Y_0\|-\|\Phi_{\cal A}(Y_0)\|}{\|Y_0\|}\geq
\frac{(j-k)^2-4\sqrt{n}m-2n}{2m^2}.$$

For each positive integer $p$, we replace $m$ with $pm$. Note that
$$Y_0=\sum\limits_{s_1,s_{2},\cdots,s_n}\sum\limits_{s'_{1},s'_{2},\cdots,s'_n}F^{pm}_{s_1,s_2,\cdots,s_n}Y_0F^{mp}_{s_1',s_2'\cdots,s_n'}\neq
0,$$ so there exist $s_1,s_{2},\cdots,s_n$ and
$s'_{1},s'_{2},\cdots,s'_n$ such that
$$Y=F^{pm}_{s_1,\cdots,s_n}Y_0F^{pm}_{s_1',\cdots,s_n'}\neq 0.$$ Thus, it is easily to prove that
$\frac{k_i}{m}\leq\frac{s_i}{pm}\leq \frac{k_i+1}{m}$ and
$\frac{k'_i}{m}\leq\frac{s'_i}{pm}\leq \frac{k'_i+1}{m}$. Note that
$k_1=k, k_1'=j$ and $|\frac{j-k}{m}|\geq \frac{2}{m}$, we have
$$\|\frac{s_1-s'_1}{pm}\|\geq  \|\frac{k_1+1-k'_1}{m}\| \geq 1/m,$$
thus $$\|s_1-s'_1\|\geq p.$$ By the similar analysis methods as (5),
we get
\begin{equation}
2(pm)^2(1-\sum\limits_{i=1}^{n}\frac{s_is'_i}{(pm)^2}-\sum\limits_{i=1}^{n}\frac{s_i+s'_i}{(pm)^2}-\frac{n}{(pm)^2})
\geq p^2-4\sqrt{n}mp-2n.\\
\end{equation} On the other hand, we also have
$$\|Y\|-\|\sum\limits_{i=1}^{n}E_iYE_i\|\geq
\|Y\|-\sum\limits_{i=1}^{n}\|E_iYE_i\|$$ $$\geq
[1-(\sum\limits_{i=1}^{n}\frac{k_ik'_i}{m^2}+\sum\limits_{i=1}^{n}\frac{k_i+k'_i}{{(pm)}^2}+\frac{n}{{(pm)}^2})]\|Y\|.$$

Let $P=F^{pm}_{s_1,s_2,\cdots,s_n}P_0$ and
$Q=Q_0F^{pm}_{s_1',s_2',\cdots,s_n'}$. Then it is clear that $P,Q\in
\mathcal{A}'$, $PQ=0$, $Y=PXQ\not=0$, and
$$\frac{\|Y\|-\|\Phi_{\cal A}(Y)\|}{\|Y\|}\geq
\frac{p^2-4\sqrt{n}m-2n}{2{(pm)}^2}.$$ The lemma is proved.

\vskip 0.1 in

It follows from the proof of Lemma 2.4 that we have the following
important conclusion:

\vskip 0.1 in

{\bf Corollary 2.1.} Let $1\leq n < \infty$, ${\cal
A}=\{E_i\}_{i=1}^{n}\subseteq {\cal E}(H)$ be commutative and
$\sum_{i=1}^{n}E_i^2\leq I$. If $X\in B(H)$ and there exist integers
$m,k,j$ with $|k-j|\geq 2$ such that
$$P^{E_1}(\frac{k}{m},\frac{k+1}{m}]XP^{E_1}(\frac{j}{m},\frac{j+1}{m}]\neq
0,$$ then for each positive integer $p$, there exist projection
operators $P,Q\in \mathcal{A}'$, $PQ=0$, $Y=PXQ\not=0$, and
$$\frac{\|Y\|-\|\Phi_{\cal A}(Y)\|}{\|Y\|}\geq
\frac{p^2-4\sqrt{n}mp-2n}{2(pm)^2}.$$

\vskip 0.2 in

\noindent {\bf 3. Main results and proofs}

\vskip 0.2 in

Let ${\cal A}=\{E_i\}_{i=1}^{n}\subseteq {\cal E}(H)$ and
$\Phi_{{\cal A}}$ be the L\"{u}ders quantum operation which is
decided by $\cal A$. It is easy to prove that $||\Phi_{{\cal
A}}||=||\sum\limits_{i=1}^n E_i^2||$ ([9]). Now, we denote
$B(H)^{\Phi_{{\cal A}}}$ to be the fixed point set of $\Phi_{{\cal
A}}$, and ${{\cal A}}'$ to be the commutant of ${{\cal A}}$, that
is, $B(H)^{\Phi_{{\cal A}}}=\{B\in B(H)\mid \Phi_{{\cal A}}(B)=B\}$,
${{\cal A}}'=\{B\in {\cal B}(H)\mid BE_i=E_i B, 1\leq i\leq n\}$. It
is clear that if $\sum_{i=1}^{n}E_i^2=I$ in strong operator
topology, then ${{\cal A}}'\subseteq B(H)^{\Phi_{{\cal A}}}$.

\vskip 0.1 in

{\bf Theorem 3.1}. Let $1\leq n\leq\infty$, ${\cal
A}=\{E_i\}_{i=1}^{n}\subseteq {\cal E}(H)$ be commutative and
$\sum_{i=1}^{n}E_i^2=I$ in strong operator topology. Then
$${\cal B}(H)^{\Phi_{\mathcal{A}}}=\{B\in {\cal
B}(H)|\Phi_{\mathcal{A}}(B)= \sum\limits_{i=1}^{n}E_iB
E_i=B\}=\mathcal {A}'.$$

{\bf Proof.} Since $\mathcal {A}'\subseteq {\cal
B}(H)^{\Phi_{\mathcal{A}}}$, in order to prove the converse
containing relation, we suppose that $B\in{\cal B}
(H)^{\Phi_{\mathcal{A}}}\setminus \mathcal {A}'$. Without of losing
generality, we can suppose that $B$ is not commutative with $E_1$.
By Lemma 2.3, there is a triple integer set $\{m,j,k\}$ such that
$|k-j|\geq 2$ and
$P^{E_1}(\frac{k}{m},\frac{k+1}{m}]BP^{E_1}(\frac{j}{m},\frac{j+1}{m}]\neq
0$.

For each positive integer $q\leq n$, let $F_q=\sum_{i=1}^qE_i^2$ and
$\Phi_q: {\cal B}(H)\rightarrow {\cal B}(H)$ be defined by
$\Phi_q(A)=\sum_{i=1}^{q}E_iAE_i$. Then $F_q\rightarrow I$ in strong
operator topology and $\Phi_q$ is a completely positive map. If
denote $P_q=P^{F_q}((1-\frac{1}{4m^2}, 1])$, then $P_q\rightarrow I$
in strong operator topology (see [18, $P_{248}]$). Now we show that
$P_qP^{E_1}(\frac{k}{m},\frac{k+1}{m}]BP^{E_1}(\frac{j}{m},\frac{j+1}{m}]P_q=0$.
In fact, if not, note that
$$P^{E_1}(\frac{k}{m},\frac{k+1}{m}]P_qP^{E_1}(\frac{k}{m},\frac{k+1}{m}]BP^{E_1}(\frac{j}{m},\frac{j+1}{m}]P_qP^{E_1}(\frac{j}{m},\frac{j+1}{m}]$$$$=P_qP^{E_1}(\frac{k}{m},\frac{k+1}{m}]BP^{E_1}(\frac{j}{m},\frac{j+1}{m}]P_q\neq
0,$$ so, by Corollary 2.1 that for each positive integer $p$, there
exist projection operators $P$ and $Q$, $P,Q\in \mathcal{A}'$,
$PQ=0$, such that
$$Y=PP_qP^{E_1}(\frac{k}{m},\frac{k+1}{m}]BP^{E_1}(\frac{j}{m},\frac{j+1}{m}]P_qQ$$$$=P_qPP^{E_1}(\frac{k}{m},\frac{k+1}{m}]BP^{E_1}(\frac{j}{m},\frac{j+1}{m}]QP_q\not
=0,$$ and
$$\frac{\|Y\|-\|\Phi_q(Y)\|}{\|Y\|}\geq
\frac{p^2-4\sqrt{q}mp-2q}{2(pm)^2}.$$

Since $$\frac{p^2-4\sqrt{q}mp-2q}{2(pm)^2}\rightarrow
\frac{1}{2m^2}$$ as $p \rightarrow \infty.$ So we can choose $Y$
such that $$\frac{\|Y\|-\|\Phi_q(Y)\|}{\|Y\|}\geq \frac{3}{8m^2}.$$
Note that $P_qE_i=E_iP_q$ and $P_qY=YP_q$ for each $1\leq i\leq n$,
${\cal A}_1=\{P_qE_i\}_{i=q+1}^n$ decides a {\it L\"{u}ders}
operation $\Phi_{{\cal A}_1}$, and
$$||\Phi_{{\cal
A}_1}||=\|\sum\limits_{i=q+1}^{n}P_qE_i^2P_q\|=\|P_q(\sum\limits_{i=q+1}^{n}E_i^2)P_q\|=\|P_q(I-\sum\limits_{i=1}^{q}E_i^2)P_q\|\leq
\frac{1}{4m^2},$$ so we have

\begin{equation}
\begin{array}{rcl}
\|\Phi_{\cal A}(Y)\|&=&\|\Phi_q(Y)+\sum\limits_{i=q+1}^{n}E_iYE_i\|\\&=&\|\Phi_q(Y)+\sum\limits_{i=q+1}^{n}E_iP_qYP_qE_i\|\\
&\leq&\|\Phi_q(Y)\|+\|\sum\limits_{i=q+1}^{n}P_qE_iYE_iP_q\|\\
&=&\|\Phi_q(Y)\|+||\Phi_{{\cal
A}_1}(Y)||\\&\leq&(1-\frac{3}{8m^2})\|Y\|+\frac{1}{4m^2}\|Y\|\\
&=& (1-\frac{1}{8m^2})\|Y\|.
\end{array}
\end{equation}

On the other hand, we show that
$Y=P_qPP^{E_1}(\frac{k}{m},\frac{k+1}{m}]BP^{E_1}(\frac{j}{m},\frac{j+1}{m}]QP_q
\in {\cal B}(H)^{\Phi_{{\cal A}}}$. In fact, note that $\{P_q, P,
P^{E_1}(\frac{k}{m},\frac{k+1}{m}],
P^{E_1}(\frac{j}{m},\frac{j+1}{m}], Q\}\subseteq \mathcal{A}'$ and
$\Phi_{{\cal A}}(B)=B$, so we have $$\Phi_{{\cal A}}(Y)=
\sum_{i=1}^{n}E_iYE_i=\sum_{i=1}^{n}E_iP_qPP^{E_1}(\frac{k}{m},\frac{k+1}{m}]BP^{E_1}(\frac{j}{m},\frac{j+1}{m}]QP_qE_i$$$$=P_qPP^{E_1}(\frac{k}{m},\frac{k+1}{m}](\sum_{i=1}^{n}E_iBE_i)P^{E_1}(\frac{j}{m},\frac{j+1}{m}]QP_q$$$$=
P_qPP^{E_1}(\frac{k}{m},\frac{k+1}{m}]\Phi_{{\cal
A}}(B)P^{E_1}(\frac{j}{m},\frac{j+1}{m}]QP_q$$$$=P_qPP^{E_1}(\frac{k}{m},\frac{k+1}{m}]BP^{E_1}(\frac{j}{m},\frac{j+1}{m}]QP_q=
Y.$$ This contradicts (7) and so
$P_qP^{E_1}(\frac{k}{m},\frac{k+1}{m}]BP^{E_1}(\frac{j}{m},\frac{j+1}{m}]P_q=0$.
Note that
$$P^{E_1}(\frac{k}{m},\frac{k+1}{m}]BP^{E_1}(\frac{j}{m},\frac{j+1}{m}]=\lim\limits_{q\rightarrow
\infty}P_qP^{E_1}(\frac{k}{m},\frac{k+1}{m}]BP^{E_1}(\frac{j}{m},\frac{j+1}{m}]P_q$$
in strong
operator topology ([17]), so $$P^{E_1}(\frac{k}{m},\frac{k+1}{m}]BP^{E_1}(\frac{j}{m},\frac{j+1}{m}]=0.$$ This contradicts $P^{E_1}(\frac{k}{m},\frac{k+1}{m}]BP^{E_1}(\frac{j}{m},\frac{j+1}{m}]\neq 0$. So $B\in \mathcal {A}'$.\\

\vskip0.1in

{\bf Theorem 3.2.} Let $1\leq n\leq\infty$, ${\cal
A}=\{E_i\}_{i=1}^{n}\subseteq {\cal E}(H)$ be commutative and
$F=\sum_{i=1}^{n}E_i^2<I$. If $P=P^F\{1\}$, where $P^F$ is the
spectral measure of $F$, then
$${\cal B}(H)^{\Phi_{\mathcal{A}}}=\{B\in {\cal
B}(H)|\Phi_{\mathcal{A}}(B)= \sum\limits_{i=1}^{n}E_iB
E_i=B\}=P\mathcal {A}'.$$

{\bf Proof.} Firstly, by the spectral representation theorem ([17])
we have $PF=FP=P$. Let $B\in {\cal B}(H)^{\Phi_{\mathcal{A}}}$. Then
as the analysis of Theorem 3.1, we have $B\in \mathcal {A}'$. Let
$Q=I-P$ and $Q_k=P^F(0,1-\frac{1}{k}]$. Then $Q_k\rightarrow Q$ in
strong operator topology and $Q_k\in \mathcal {A}'$, so $Q_kB\in
{\cal B}(H)^{\Phi_{\mathcal{A}}}$. Let $\Phi_k$ be the completely
positive map which is decided by $\{E_iQ_k\}_{i=1}^{n}$. Then
$\|\Phi_k\|\leq 1-\frac{1}{k}$. Note that $B, Q_k\in \mathcal {A}'$
and $Q_k^2=Q_k$, thus we have
$\|Q_kB\|=\|\Phi_{\mathcal{A}}(Q_kB)\|=\|\Phi_k(Q_kB)\|\leq
(1-\frac{1}{k})\|Q_kB\|$, so $Q_kB=0$. Note that
$QB=\lim\limits_{k\rightarrow \infty}Q_kB$ in strong operator
topology, so $QB=0$, that is, $(I-P)B=0$, i.e., $B=PB$, this showed
that ${\cal B}(H)^{\Phi_{\mathcal{A}}}\subseteq P\mathcal {A}'$. If
$B\in P\mathcal {A}'$, note that $P\in \mathcal {A}'$, so $PB=BP=B$.
Moreover, $\Phi_{{\cal A}}(B)=BF=PBF=BPF=BP=B$, that is, $B\in {\cal
B}(H)^{\Phi_{\mathcal{A}}}$, thus we have $P\mathcal {A}'\subseteq
{\cal B}(H)^{\Phi_{\mathcal{A}}}$ and the theorem is proved.

\vskip0.2in

{\bf Acknowledgement.} The authors wish to express their thanks to
the referees for their valuable comments and suggestions. This
project is supported by Zhejiang Innovation Program for Graduates
(YK2009002) and Natural Science Foundations of China (10771191 and
10471124) and Natural Science Foundation of Zhejiang Province of
China (Y6090105).

\vskip0.2in

\centerline{\bf References}

\vskip0.2in

\noindent [1]. Foulis, D. J., Bennett, M. K. Effect algebras and
unsharp quantum logics. Found. Phys. 24(10), 1331-1352 (1994)

\noindent [2] Ludwig, G. Foundations of Quantum Mechanics (I-II),
Springer, New York, 1983.

\noindent [3] Ludwig, G. An Axiomatic Basis for Quantum Mechanics
(II), Springer, New York, 1986.

\noindent [4] Kraus, K. Effects and Operations, Springer, New York,
1983.

\noindent [5] Davies, E. B. Quantum Theory of Open Systems, Academic
Press, London, 1976.

\noindent [6]. Busch, P, Grabowski, M. and Lahti, P. J. Operational
Quantum Physics, Springer-Verlag, Beijing Word Publishing
Corporation (1999)

\noindent [7]. L\"{u}ders, G. \"{U}ber die Zustands\"{a}nderung
durch den Messprozess, Ann. Physik (6). 8 (1), 322-328 (1951)

\noindent [8]. Busch, P and Singh, J. L\"{u}ders theorem for unsharp
quantum measurements. Phys. Letter A. 249 (1-2), 10-12 (1998)

\noindent [9]. Arias, A., Gheondea, A. and Gudder, S. Fixed points
of quantum operations. J. Math. Phys. 43(12), 5872-5881 (2002)

\noindent [10]. Kribs, D. W.: Quantum channels, wavelets, dilations,
and representations of ${\cal O}_n$. Proc. Edinb. Math. Soc. 46 (2),
421-433 (2003)

\noindent [11]. Liu Weihua, Wu Junde. On fixed points of L\"{u}ders
operation. J. Math. Phys. 50(10), 103531-103532(2009)

\noindent [12]. Nagy, G. On spectra of L\"{u}ders operations. J.
Math. Phys. 49 (2), 022110-022117 (2008)

\noindent [13]. Holbrook, J. A., Kribs, D. W., Laflamme, R. and
Poulin, D. Noiseless subsystems for collective rotation channels in
quantum information theory. Inter. Equ. Oper. Theory. 51 (2),
215-234 (2005)

\noindent [14]. Choi, M. D, Kribs D. W. Method to Find Quantum
Noiseless Subsystems. Phys. Rev. Lett. 96 (5), 050501-050504 (2006)

\noindent [15]. Blume-Kohout R., Ng, H. K, Poulin, D. and Viola, L.
Characterizing the Structure of Preserved Information in Quantum
Processes. Phys. Rev. Lett. 100 (3), 030501-030504 (2008)

\noindent [16]. Choi, M. D., Johnston, N., Kribs, D. W. The
multiplicative domain in quantum error correction. J. Phys. A: Math.
Theor. 42 (24), 245303-245317 (2009)

\noindent [17]. Kadison, R. V. and Ringrose, J. R.: Fundamentals of
the Theory of Operator algebra I, II, Springer, New York (1983)

\noindent [18]. Riesz, F. and SZ.-Nagy, B.: Functional Analysis,
Science Press of China, Beijing (1981)

\end{document}